\documentclass[12pt]{article}
\usepackage{amsmath,amssymb,amscd,latexsym}
\usepackage[all]{xy}

\newcommand{\C}{{\mathbb{C}}}
\newcommand{\F}{{\mathbb{F}}}

\newcommand{\Q}{{\mathbb{Q}}}
\newcommand{\oQ}{\overline{\Q}}

\newcommand{\R}{{\mathbb{R}}}
\newcommand{\Z}{{\mathbb{Z}}}

\newcommand{\hL}{\hat{L}}

\newcommand{\ddet}{\mathrm{det}}

\newcommand{\et}{\mathrm{\acute{e}t}}

\newcommand{\id}{\mathrm{id}}
\newcommand{\ord}{\mathrm{ord}}
\newcommand{\rk}{\mathrm{rk}}
\newcommand{\sing}{\mathrm{sing}}

\newcommand{\spec}{\mathrm{spec}\,}

\newcommand{\Ker}{\mathrm{Ker}\,}

\newcommand{\tL}{\tilde{L}}
\newcommand{\ty}{\tilde{y}}

\newcommand{\RRe}{\mathrm{Re}\,}

\newcommand{\tr}{\mathrm{tr}}

\newcommand{\Ch}{{\mathcal C}}

\newcommand{\Fh}{{\mathcal F}}

\newcommand{\eX}{{\mathfrak X}}

\newcommand{\oeX}{\overline{\eX}}

\newcommand{\hzeta}{\hat{\zeta}}

\newcommand{\oF}{\overline{\mathbb{F}}}

\newcommand{\tc}{\tilde{c}}

\newcommand{\tx}{\tilde{x}}
\newcommand{\orho}{\overline{\rho}}
\newcommand{\olambda}{\overline{\lambda}}
\newcommand{\tei}{\, | \,}
\newcommand{\verk}{\mbox{\scriptsize $\,\circ\,$}}

\newtheorem{theorem}{Theorem}

\newtheorem{example}[theorem]{Example}
\newtheorem{remark}[theorem]{Remark}

\newtheorem{conjecture}[theorem]{Conjecture}
\begin{document}
\title{The Hilbert--Polya  strategy and height pairings}
\author{C. Deninger}
\date{\ }
\maketitle

\begin{abstract} Previously we gave a conjectural cohomological argument for the validity of the Riemann hypotheses for Hasse--Weil zeta functions. In the present note we sketch how the same cohomological formalism would imply the conjectured positivity properties of the height pairings of homologically trivial cycles.
\end{abstract}

\section{Introduction}
\label{sec:1}

In earlier work which will be recalled below we have introduced a conjectural cohomology theory which would allow a natural proof of the Riemann hypotheses for the zeta- and $L$-functions occuring in arithmetic geometry: The argument used a scalar product on cohomology involving the cup product and a Hodge $*$-operator. In this note we show that these ingredients would also help one to prove that the height pairing on homologically trivial cycles has the non-degeneracy and positivity properties conjectured by Beilinson in \cite{Be} and Bloch in \cite{Bl}. By the work of N\'eron and Tate, the positivity of the height pairing is known for divisors at least. Hence our remark may be interpreted as support for the Hilbert--Polya strategy of proving the Riemann hypotheses via positivity.

\section{The conjectural cohomological formalism and the height pairing}
\label{sec:2}

For a projective scheme $\eX / \Z$ consider the Hasse--Weil zeta function
\[
\zeta_{\eX} (s) = \prod_{x \in |\eX|} (1 - (Nx)^{-s})^{-1} \; .
\]
The product runs over the set $|\eX|$ of closed points of $\eX$ and $Nx$ is the number of elements in the residue field of $x$. The product converges absolutely and locally uniformly for $\RRe s > d = \dim \eX$ and hence defines an analytic function in that half plane.
For $\eX = \spec \Z$ we get the Riemann zeta function for example. 

The completed zeta function $\hat{\zeta}_{\eX} (s)$ is defined by taking the product of $\zeta_{\eX} (s)$ with a finite number of factors involving the Euler Gamma function. The precise receipe for which we refer to \cite{S2} depends on the real Hodge structure of $\eX \otimes \R$. For example we have:
\[
\hzeta (s) = \zeta (s) \pi^{-s/2} \Gamma (s/2) \; .
\]
If $\eX$ is defined over a finite field, so that $\eX \otimes \R$ is empty we have $\hzeta_{\eX} (s) = \zeta_{\eX} (s)$.

A lot of work by many mathematicians has led to the following conjecture where we now assume in addition that $\eX$ is regular

\begin{conjecture}
  \label{t1}
a) $\zeta_{\eX} (s)$ has a meromorphic continuation to all of $\C$.\\
b) There is a functional equation of the form:
\[
\hzeta_{\eX} (d-s) = A B^s \hzeta_{\eX} (s)
\]
for certain constants $A$ and $B > 0$ c.f. \cite{S2}.\\
c) (Riemann hypotheses) The zeroes (poles) of $\hzeta_{\eX} (s)$ lie on the lines $\RRe s = \frac{i}{2}$ for $0 \le i \le 2d$ odd (even).
\end{conjecture}

If $\eX$ is an $\F_p$-scheme the conjecture has been proved by Grothendieck, Deligne and many others c.f. \cite{D}. The proof uses the $l$-adic cohomology theory of Grothendieck. If $\eX / \Z$ is flat, conjectures a) and b) are known in several cases by first relating $\zeta_{\eX} (s)$ via trace formulas to automorphic $L$-functions and then using their theory. Conjecture c) on the other hand is not known for a single flat $\Z$-scheme $\eX$.

Clearly it would be very desirable to have a cohomology theory which would serve the same purposes also for flat schemes $\eX / \Z$ as the $l$-adic cohomology does for $\F_p$-schemes $\eX$. In this regard we have the following conjecture whose compatibility with known or expected properties of $\zeta_{\eX} (s)$ has been checked in \cite{D1}, \cite{D2}.

\begin{conjecture}
  \label{t2}
On the category of separated schemes of finite type $\eX / \Z$ there exists a cohomology theory $H^i (\oeX)$ of complex Fr\'echet spaces with $\R$-action $\phi^{t*}$ such that for the infinitesimal operator
\[
\theta = \lim_{t\to 0} \frac{1}{t} (\phi^{t*} - \id)
\]
we have
\[
\hzeta_{\eX} (s) = \prod^{2d}_{\nu=0} \ddet_{\infty} \Big( \frac{1}{2\pi} (s\, \id - \theta) \tei H^{\nu} (\oeX) \Big)^{(-1)^{\nu+1}} \; .
\]
\end{conjecture}

Here the spectrum of $\theta$ on $H^{\nu} (\oeX)$ is supposed to consist of eigenvalues only and $\det_{\infty}$ is the zeta-regularized determinant of the endomorphism $\varphi = \frac{1}{2\pi} (s \, \id -\theta)$ acting on $H^{\nu} (\oeX)$. The eigenvalues of $\varphi$ are the numbers $\lambda = \frac{1}{2\pi} (s - \rho)$ where $\rho$ runs over the eigenvalues of $\theta$. Let
\[
\eta (z) = \sum_{\lambda \neq 0} \lambda^{-z} = \sum_{\rho \neq s} \Big[ \frac{1}{2\pi} (s- \rho) \Big]^{-z}
\]
be the spectral zeta function of $\varphi$. Here the eigenvalues occur with their algebraic multiplicities, which are supposed to be finite. Then $\det_{\infty} (\varphi)$ is defined if $\eta$ converges in some right half plane with an analytic continuation to $\RRe z > - \varepsilon$ for some $\varepsilon > 0$. In this case one sets $\det_{\infty} (\varphi) = \exp (- \eta' (0))$ if $0 \notin \spec \varphi$ i.e. if $s \notin \spec \theta$ and $\det_{\infty} (\varphi) = 0$ otherwise. 

As evidence for conjecture \ref{t2} we note that it is true if $\eX$ is an $\F_p$-scheme. The construction of $H^{\nu} (\oeX)$ in that case given in \cite{D1} \S\,3 is unsatisfactory however because it is not geometric but instead based on the already existing $l$-adic cohomology groups $H^{\nu}_{\et} (\eX \otimes \oF_p , \Q_l)$. Moreover an embedding of $\Q_l$ into $\C$ has to be fixed which looks very artificial. But in fact something similar happens in a familiar situation: because of the comparison isomorphisms
\[
H^{\nu}_{dR} (X / \C) = H^{\nu}_{\sing} (X (\C) , \Q) \otimes \C \quad \mbox{and} \quad H^{\nu}_{\et} (X \otimes \oQ , \Q_l) = H^{\nu}_{\sing} (X (\C) , \Q) \otimes \Q_l
\]
every embedding $\sigma : \Q_l \hookrightarrow \C$ induces an isomorphism
\[
H^{\nu}_{dR} (X / \C) = H^{\nu}_{\et} (X \otimes \oQ , \Q_l) \otimes_{\Q_l, \sigma} \C \; .
\]
Thus, if de~Rham cohomology were not known yet, it could be obtained from $l$-adic cohomology after choosing an embedding of $\Q_l$ into $\C$.

One of the consequences of conjecture \ref{t2} is that $\hzeta_{\eX} (s)$ should be zeta-regularizable i.e. obtained from its divisor of zeros and poles by a canonical procedure involving zeta-regularization. For the Riemann zeta function this was proved in \cite{ScSo}. Later Illies \cite{I}, see also \cite{D3} \S\,5, showed that it follows in general from some of the other expected analytic properties of $\zeta_{\eX} (s)$.

The Riemann hypotheses for $\hzeta_{\eX} (s)$ would correspond naturally to the relation $\RRe \rho = \nu / 2$ for the eigenvalues $\rho$ of $\theta$ on $H^{\nu} (\oeX)$.

Let us now give a conjectural argument in the spirit of \cite{S1} for the Riemann hypotheses. Firstly we should have an $\R$-equivariant trace homomorphism
\[
\tr : H^{2d} (\oeX) \to \C (-d)
\]
where $t \in \R$ acts on $H^{2d} (\oeX)$ by $\phi^{t*}$ and on $\C (-d) := \C$ by multiplication with $e^{dt}$. This is reasonable since $\hzeta_{\eX} (s)$ is expected to have a pole at $s = d$, i.e. $\theta$ should have the eigenvalue $\rho = d$ on $H^{2d} (\oeX)$.

Next we expect an $\R$-equivariant cup-product pairing:
\[
H^{\nu} (\oeX) \times H^{2d-\nu} (\oeX) \xrightarrow{\cup} H^{2d} (\oeX) \xrightarrow{\tr} \C (-d) \; .
\]
The induced pairing between the $\rho$-eigenspace of $\theta$ on $H^{\nu} (\oeX)$ and the $d-\rho$ eigenspace of $\theta$ on $H^{2d-\nu} (\oeX)$ should be perfect for all $\rho \in \C$. In particular $\rho$ would be an eigenvalue of $\theta$ on $H^{\nu} (\oeX)$ if and only if $d-\rho$ were an eigenvalue of $\theta$ on $H^{2d-\nu} (\oeX)$. This is compatible with the conjectured functional equation of $\hzeta_{\eX} (s)$. 

Finally we expect a $\C$-antilinear Hodge $*$-isomorphism
\[
* : H^{\nu} (\oeX) \xrightarrow{\sim} H^{2d-\nu} (\oeX)
\]
such that the hermitian form on $H^{\nu} (\oeX)$
\[
(f , f') = \tr (f \cup * (f'))
\]
is positive definite i.e. a scalar product. The $*$-operator and the $\R$-action should be related by the formula
\[
\phi^{t*} \verk * = e^{t (d-\nu)} * \verk \phi^{t*} \quad \mbox{on} \; H^{\nu} (\oeX) \; .
\]
Equivalently we would have
\[
\theta \verk * = * \verk (d-\nu + \theta) \quad \mbox{on} \; H^{\nu} (\oeX) \; .
\]
In the foliation analogue of \cite{D4} this relation holds on foliation cohomology if there is a flow $\phi^t$ inducing $\phi^{t*}$ which is conformal with factor $e^t$ for a leafwise metric. From that perspective the relation is quite natural therefore.

The above formulas would imply the Riemann hypotheses as follows: For $f, f'$ in $H^i (\oeX)$, since $\theta$ is a derivation with respect to cup product, we would have:
\begin{eqnarray*}
  d (f,f') & = & \tr (\theta (f \cup * f') = \tr (\theta f \cup * f') + \tr (f \cup \theta * f') \\
 & = & (\theta f , f') + (d-\nu) (f, f') + (f , \theta f') \; .
\end{eqnarray*}
The relation
\[
\nu (f, f') = (\theta f, f') + (f , \theta f')
\]
would then show that $\theta - \nu/2$ is skew symmetric on $H^{\nu} (\oeX)$ so that its spectrum would be purely imaginary.

For the comparison with the conjectured positivity properties of the height pairing we now look at a version of this argument for a slightly different cohomology that should be relevant for an understanding of motivic $L$-functions.

Given a smooth projective scheme $X / \Q$ of equidimension $\delta$ one defines the $w$-th $L$-function $L (H^w (X), s)$ and its completed version $\hL (H^w (X), s)$ as in \cite{S2}. If $X = E$ is an elliptic curve over $\Q$ for example we have
\[
L (H^1 (E) , s) = L (E,s)
\]
the $L$-function of $E$ and
\[
\hL (H^1 (E) , s) = L (E,s) (2 \pi)^{-s} \Gamma (s) \; .
\]
Granting the analytic continuation, the Riemann hypotheses for these $L$-functions asserts that the zeroes of $\hL (H^w (X) , s)$ should satisfy $\RRe s = \frac{w+1}{2}$ and the poles $\RRe s = \frac{w}{2}$ or $\RRe s = \frac{w}{2} + 1$. Concerning the poles there is the more precise Artin conjecture, which asserts that $\hL (H^w (X) , s)$ can have a pole if and only if $w$ is even and the motive $H^w (X)$ contains the Tate motive $\Q (-w/2)$ as a direct factor. In that case $\hzeta (s - \frac{w}{2})$ is a factor and $\hL (H^w (X) , s)$ would have poles precisely at the points $s = \frac{w}{2}$ and $s = \frac{w}{2} + 1$. 

The conjectural cohomological formalism of \cite{D1} predicts a formula of the form
\[
\hL (H^w (X) , s) = \prod^2_{\nu=0} \ddet_{\infty} \Big( \frac{1}{2\pi} (s \, \id - \theta) \tei H^{\nu} (\overline{\spec \Z} , \Fh^{w}_X) \Big)^{(-1)^{\nu+1}}
\]
where the sheaf $\Fh^w_X$ should come about as follows: Choose a smooth projective model $\pi : \eX_U \to U$ of $X$ over some non-empty open subscheme $U$ of $\spec \Z$ and set $\Fh^w_X = j_* R^w \pi_* \Ch$ where $j : U \hookrightarrow \overline{\spec \Z}$ is the open immersion and $\Ch$ the structural sheaf in the conjectured site. The sheaf $\Fh^w_X$ should depend only on $X$ and not on the choices of $U$ and the model $\eX_U$ over $U$. Everything  should be equipped with an $\R$-action $\phi^{t*}$ and $\theta$ denotes the infinitesimal generator. 

Cup product and trace maps should give pairings
\[
\Fh^w_X \times \Fh^{2\delta - w}_X \xrightarrow{\cup} j_* R^{2\delta} \pi_* \Ch \xrightarrow{\tr} j_* \Ch (-\delta)
\]
and
\[
H^{\nu} (\overline{\spec \Z} , \Fh^w_X) \times H^{2-\nu} (\overline{\spec \Z} , \Fh^{2\delta -w}_X) \to H^2 (\overline{\spec \Z} , j_* \Ch) (-\delta) \; .
\]
In our formalism the group $H^{\nu} (\overline{\spec \Z})$ from above for $\eX = \spec \Z$ is actually an abbreviation of $H^{\nu} (\overline{\spec \Z} , j_* \Ch)$. So together with the trace map on $H^2 (\overline{\spec \Z})$ we would obtain a pairing:
\[
H^{\nu} (\overline{\spec \Z} , \Fh^w_X) \times H^{2-\nu} (\overline{\spec \Z} , \Fh^{2\delta -w}_X) \xrightarrow{\cup} H^2 (\overline{\spec \Z}) (-\delta) \xrightarrow{\tr} \C (-d)
\]
where $d = \delta+1$ is the dimension of $\eX$.

Now assume that there is a $\C$-antilinear $*$-operator
\[
* : H^{\nu} (\overline{\spec \Z} , \Fh^w_X) \xrightarrow{\sim} H^{2-\nu} (\overline{\spec \Z} , \Fh^{2 \delta - w}_X)
\]
such that setting
\[
(f,f') = \tr (f \cup *f')
\]
defines a scalar product on $H^{\nu} (\overline{\spec \Z} , \Fh^w_X)$. If on this group we had
\[
\phi^{t*} \verk * = e^{t (d-\nu-w)} * \verk \phi^{t*}
\]
i.e.
\[
\theta \verk * = * \verk (d-\nu-w + \theta)
\]
then arguing as before we would get
\[
(\nu+w) (f, f') = (\theta f ,f ') + (f , \theta f') \; .
\]
Hence $\theta - \frac{\nu+w}{2}$ would be skew symmetric on $H^{\nu} (\overline{\spec \Z} , \Fh^w_X)$ and thus the Riemann hypotheses for $\hL (H^w (X) , s)$ would follow.

Note also that one expects a functional equation of the form
\[
\hL (H^w (X) , s) = ab^s \hL (H^w (X) , w+1-s)
\]
for certain constants $a$ and $b > 0$. In particular $\rho$ would be a zero (or pole) of $L (H^w (X) , s)$ if and only if $w+1-\rho$ is a zero (or pole). This is compatible with the following argument: If $\rho$ is an eigenvalue of $\theta$ on $H^{\nu} (\overline{\spec \Z} , \Fh^w_X)$ then by Poincar\'e duality $d-\rho$ should be an eigenvalue of $\theta$ on $H^{2-\nu} (\overline{\spec \Z} , \Fh^{2\delta -w}_X)$. By the $*$-isomorphism it would follow that $\nu +w-\orho$ is an eigenvalue of $\theta$ on $H^{\nu} (\overline{\spec \Z} , \Fh^w_X)$ as well. This is compatible with the functional equation since $\lambda$ is a zero of $L(H^w (X) , s)$ if and only if $\olambda$ is a zero.

Now let us discuss the compatibility of this picture with the height pairing on algebraic cycles. Following Beilinson's notation in a related setting we write:
\[
H^p_{!*} (X) = H^1 (\overline{\spec \Z} , \Fh^{p-1}_X) \; .
\]
In this notation the above statements read as follows. We have a pairing given by cup product and the trace map:
\[
\langle , \rangle : H^p_{!*} (X) \times H^{2d-p}_{!*} (X) \to \C (-d)
\]
which should be perfect on $\theta$-eigenspaces. The $*$-operator would give an antilinear isomorphism:
\[
* : H^p_{!*} (X) \xrightarrow{\sim} H^{2d-p}_{!*} (X) 
\]
such that $\theta \verk * = * \verk (d-p+\theta)$ on $H^p_{!*} (X)$. Moreover the ``Hodge inner product'' $(f,f') = \langle f, *f' \rangle$ would be positive definite on $H^p_{!*} (X)$. 

The generalized Birch, Swinnerton-Dyer conjecture asserts the following equality
\[
\ord_{s=i} L (H^{2i-1} (X) , s) = \dim CH^i (X)^0 \; .
\]
Here
\[
CH^i (X)^0 = \Ker (CH^i (X) \otimes \Q \to H^{2i} (X (\C) , \Q))
\]
denotes the $\Q$-vector space of codimension $i$ cycles on $X$ up to rational equivalence whose cohomology classes are trivial.

Assuming the above formalism we have:
\begin{eqnarray*}
  \ord_{s=i} L (H^{2i-1} (X) , s) & = & \ord_{s=i} \hL (H^{2i-1} (X) , s) \\
& = & \dim_{\C} H^1 (\overline{\spec \Z} , \Fh^{2i-1}_X)^{\theta=i} \\
& = & \dim_{\C} H^{2i}_{!*} (X)^{\theta = i} \; .
\end{eqnarray*}
Hence the conjecture would follow from an Abel--Jacobi isomorphism which we write as $c \mapsto \tc$
\[
^{\sim} : CH^i (X)^0 \otimes \C \xrightarrow{\sim} H^{2i}_{!*} (X)^{\theta = i} \; .
\]
The height pairing $\langle , \rangle_{CH}$ on cycles which are homologically equivalent to zero \cite{Be} should fit into the following diagram where $\langle f_1 , f_2 \rangle = \tr (f_1 \cup f_2)$
\[
\xymatrix{
CH^i (X)^0 \ar[d] \ar@{}[r] |{\times} & CH^{d-i} (X)^0 \ar[d] \ar[r]^-{\langle , \rangle_{CH}} & \R \ar[d] \\
H^{2i}_{!*} (X)^{\theta = i} \ar@{}[r] |{\times} & H^{2d - 2i}_{!*} (X)^{\theta = d-i} \ar[r]^-{\langle , \rangle} & \C \; .
}
\]

The following conjecture due to Bloch and Beilinson is proved in the case $i = 1$ by the work of N\'eron and Tate

\begin{conjecture}[\cite{Be} \S\,5]
  \label{t3}
a) $\langle , \rangle_{CH}$ is non-degenerate.\\
b) Let $L \in CH^1 (X)$ be the class of a hyperplane section. Then for $i \le d/2$ the map
\[
L^{d-2i} : CH^i (X)^0 \xrightarrow{\sim} CH^{d-i} (X)^0
\]
is an isomorphism. On the primitive $i$-cycles $x$ i.e. those with $i \le d/2$ and $L^{d-2i+1} (x) = 0$ the pairing
\[
(x,y)_{CH} = \langle x, L^{d-2i} y \rangle_{CH}
\]
is definite of sign $(-1)^i$.
\end{conjecture}

\begin{remark}
  In the cohomological theory a cohomology class $h$ of degree $p$ is primitive if $p \le d$ and $L^{d-p+1} (h) = 0$. This is compatible with the above if we think of an $i$-cycle in $CH^i (X)^0$ as having cohomological degree $2i$.
\end{remark}

Now assume that $L$ has a cohomology class $\tL$ which induces a hard Lefschetz isomorphism for $i \le d/2$:
\[
\tL^{d-2i} : H^{2i}_{!*} (X) \xrightarrow{\sim} H^{2d-2i}_{!*} (X) \; .
\]
If $\phi^{t*} (\tL) = e^t \tL$ i.e. $\theta (\tL) = \tL$, then this gives an isomorphism
\[
\tL^{d-2i} : H^{2i}_{!*} (X)^{\theta=i} \xrightarrow{\sim} H^{2d-2i}_{!*} (X)^{\theta = d-i} \quad \mbox{for} \; i \le d/2 \; .
\]
Of course the following diagram should commute:
\[
\xymatrix{
CH^i (X)^0 \ar[r]^{L^{d-2i}} \ar[d] & CH^{d-i} (X)^0 \ar[d] \\
H^{2i}_{!*} (X) \ar[r]^{\tL^{d-2i}} & H^{2d-2i}_{!*} (X)
}
\]
In particular primitive cycles are mapped to primitive cycles by the vertical Abel--Jacobi maps. For $x,y$ in $CH^i (X)^0$ we then get:
\[
(x,y)_{CH} = \langle x , L^{d-2i} y \rangle_{CH} = \langle \tx , \tL^{d-2i} \ty \rangle \; .
\]
Now, in classical K\"ahler theory one has the following formula for primitive $(i,i)$-forms $\psi$  in $H^{2i}$ of a complex $d$-dimensional compact K\"ahler manifold:
\[
* \psi = \frac{(-1)^i}{(d-2i)!} L^{d-2i} \psi \; .
\]
Hence for {\it primitive} $y$ we would expect:
\[
\tL^{d-2i} \ty = (-1)^i (d - 2i)! * \ty \; .
\]
This would give:
\[
(x,y)_{CH} = (-1)^i (d-2i)! (\tx , \ty) \; .
\]
Thus the positivity of $(,)$ on $H^{2i}_{!*} (X)$ which was needed in the above argument for the Riemann hpyotheses would also be responsible for the positivity properties of the height pairings.

\begin{example}
  Let us look at the simplest non-trivial case where $E / \Q$ is an elliptic curve. Then with normalizations as in \cite{Be} it is known by the work of N\'eron and Tate that the height pairing on $CH^1 (E)^0 = E (\Q) \otimes \Q$
\[
\langle , \rangle_{CH} : CH^1 (E)^0 \times CH^1 (E)^0 \to \R
\]
is negative definite. Note that here $d = 2 , i = 1$ and hence all cycles in $CH^1 (E)^0$ are primitive. The Abel--Jacobi map would be an isomorphism
\[
E (\Q) \otimes \C \xrightarrow{\sim} H^2_{!*} (E)^{\theta = \id}
\]
explaining the Birch Swinnerton-Dyer conjecture
\[
\ord_{s=1} L (E,s) = \rk \, E (\Q) \; .
\]
The relation $\theta \verk * = * \verk \theta$ on $H^2_{!*} (E)$ would show that we have $\theta = 1+S$ where $S$ is skew-symmetric with respect to the Hodge inner product. Finally the restriction of the Hodge inner product to the $1$-eigenspace of $\theta$ would be identified with the negative of the (negative-definite) N\'eron--Tate height pairing on $E (\Q) \otimes \C$.
\end{example}

There are reasons to think that if $\eX$ is flat over $\Z$ the hoped for cohomologies carry a natural real structure. In particular the cohomologies $H^{\nu} (\oeX)$ and $H^{\nu} (\overline{\spec \Z} , \Fh^w_X)$ and $H^{\nu}_{!*} (X)$ should have real structures respected by $\phi^{t*}$ and hence by $\theta$. Everything we said above is compatible with this possibility.\\

{\bf Acknowledgements.} I would like to thank Masato Wakayama for the invitation which led to the present note.

\end{document}